\documentclass{amsart}

\usepackage{amsthm}
\usepackage{graphicx}
\usepackage{amssymb}
\usepackage{psfrag}

\newtheorem{teo}{Theorem}[section]
\newtheorem{lemma}[teo]{Lemma}
\newtheorem{prop}[teo]{Proposition}
\newtheorem{cor}[teo]{Corollary}

\newtheorem{question}[teo]{Question}

\theoremstyle{definition}
\newtheorem{defi}[teo]{Definition}
\newtheorem{rem}[teo]{Remark}

\theoremstyle{remark}
\newtheorem{prof}[teo]{Proof of}

\def\mr{\mathbb{R}}
\def\mz{\mathbb{Z}}

\begin{document}
\title{Stein domains and branched shadows of $4$-manifolds}

\author[Costantino]{Francesco Costantino}
\address{Institut de Recherche Math\'ematique Avanc\'ee\\
  7, Rue Ren\'e Descartes\\
  67084 Strasbourg Cedex France}
\email{f.costantino@sns.it}

\begin{abstract}
We provide sufficient conditions assuring that a suitably decorated $2$-polyhedron can be thickened to a compact $4$-dimensional Stein domain. 
We also study a class of flat polyhedra in $4$-manifolds and find conditions assuring that they admit Stein, compact neighborhoods. We base our calculations on Turaev's shadows suitably ``smoothed"; the conditions we find are purely algebraic and combinatorial.  Applying our results, we provide examples of hyperbolic $3$-manifolds admitting ``many" positive and negative Stein fillable contact structures, and prove a $4$-dimensional analogue of Oertel's result on incompressibility of surfaces carried by branched polyhedra.
\end{abstract}

\maketitle

\tableofcontents
\section{Introduction}
The present paper is devoted to find combinatorial conditions allowing one to reconstruct compact Stein domains from branched shadows. Branched shadows, defined in \cite{Co} and \cite{Cobr}, can be considered as decorated $2$-dimensional polyhedra which can be canonically thickened to $4$-manifolds and whose singularities are smooth; equivalently, they are smooth $2$-polyhedra flatly embedded in $4$-manifolds admitting bases of neighborhoods diffeomorphic to the ambient manifolds. They are the natural combination of Turaev's shadows of $4$-manifolds (\cite{Tu}, \cite{Tu2}) and Benedetti and Petronio's branched spines of $3$-manifolds (\cite{BP}). 
In \cite{Cobr}, we showed how a branched shadow  ``carries" an almost complex structure well defined up to homotopy on its thickening and, conversely, that each almost complex structure on the thickening is carried by a suitable shadow.  We also proved that each almost complex structure on a $4$-dimensional handlebody is homotopic to a complex one.

In the present paper we prove two results relating shadows to Stein domains. The first one is a combinatorial sufficient condition on an abstract, non embedded, branched shadow assuring that its canonical thickening can be equipped with a Stein domain structure. It represents the natural translation in the world of shadows of the well-known necessary and sufficient condition provided by Gompf (\cite{Go}) and Eliashberg (\cite{El}). The following is a simplified statement of the result (Theorem \ref{teo:mainteo}):
\begin{teo}
If the canonical $2$-cochain which represents the first Chern class of the almost complex structures carried by the shadow is ``sufficiently negative", then the thickening of the shadow can be equipped with a Stein domain structure whose complex structure is homotopic to that carried by the shadow.  
\end{teo}
We stress here that the above sufficient condition is purely combinatorial. Indeed, there is a canonical $2$-cochain representing the Chern class of the complex structures carried by the shadow whose coefficients depend only on the combinatorial structure of the branched shadow. We apply the above result to the case of branched spines of $3$-manifolds, which, as showed by Benedetti and Petronio (\cite{BP}), carry a distribution of oriented $2$-planes on their thickenings.  We prove the following (see Corollary \ref{cor:branchedspines}):
\begin{teo}
The distribution of oriented $2$-planes carried by a branched spine of a $3$-manifold  whose canonical $2$-cochain is non-positive is homotopic both to a positive and to a negative tight contact structure.  
\end{teo}

In the last subsection, as an application of the above results, we exhibit families of hyperbolic $3$-manifolds admitting both positive and negative Stein fillings and prove the following:
\begin{prop}
For each integer $n\geq 1$ there exist a closed, hyperbolic $3$-manifold $Q_n$, whose volume is less than $4nVol_{Oct}$ (where $Vol_{oct}$ is the volume of the regular ideal hyperbolic octahedron) and such that $Q_n$ has at least $n-1$ positive and $n-2$ negative, Stein fillable contact structures. 
\end{prop}
We then pass to the point of view of branched polyhedra flatly embedded in complex $4$-manifolds. Roughly speaking, on each oriented smooth $2$-dimensional object in a complex manifold one can distinguish totally real, elliptic, hyperbolic positive or negative complex points. As proved by Forstneri\v c (\cite{Fo}), an oriented real surface in a complex manifold can be perturbed to one admitting a Stein tubular neighborhood provided that all the elliptic points it contains can be annihilated with (an equal number of) hyperbolic points. We study a generalization of this condition to the case of embedded branched polyhedra. To state our result, let us recall that the combinatorial structure on the polyhedron allows a canonical presentation of the cohomology of its neighborhood in the ambient manifold using as a basis the duals of the cells of the polyhedron subject to relations given by the coboundaries of one-dimension-less cells. In particular, we define two integer $2$-cochains $I^+$ and $I^-$, whose evaluation on a region is equal respectively to the total index of positive and negative complex points contained in the region. The natural translation of Forstneri\v c's condition would be that each of these $2$-cochains is cohomologous to one having only negative coefficients; unfortunately, it turns out that this is not the correct condition. Indeed, one has to deal with the behavior of the ambient complex structure near the singular set of the polyhedron. This causes that, roughly speaking, while annihilating hyperbolic and elliptic points, one has to count the number of times this annihilation ``involves" the singular set. 
In \cite{Cobr}, we proved a lemma which represents the translation in the world of shadows of the well known Harlamov-Eliashberg Annihilation Theorem for complex points over real surfaces: roughly speaking, it says that one can shift these points along a branched polyhedron, but while passing through the singular set, they ``duplicate" so that the cochains $I^+$ and $I^-$ change by a coboundary.
The resulting sufficient condition we get in Theorem \ref{teo:genfo}, can be summarized, as follows:
\begin{teo}
If there are $1$-cochains $b^+$ and $b^-$ such that $I^++\delta b^+$ and $I^-+\delta b^-$ have non-positive coefficients and $b^++b^-$ is ``sufficiently positive", then the polyhedron can be isotoped so that it admits a compact Stein tubular neighborhood.  
\end{teo}
The nice feature of the above condition is that it reduces to a purely algebraic calculation, depending only on the combinatorial structure of the branched polyhedron and on some external data induced by the embedding in the ambient manifold as, for instance, the cochains $I^{\pm}$. 

As an application of Theorem \ref{teo:genfo}, in the last subsection, we prove a result on surfaces ``carried" by a branched polyhedron i.e. surfaces lying in a regular neighborhood of the polyhedron and whose projection on it has everywhere positive differential. Our result represents the $4$-dimensional analogue of Oertel's result (\cite{Oe}) on incompressibility of surfaces carried by branched surfaces in $3$-manifolds, and states the following:
\begin{teo}
Let $(M,J)$ be a complex $4$-manifold collapsing on a branched polyhedron $P$ satisfying the hypotheses of Theorem \ref{teo:genfo} and such that $J$ is homotopic to the almost complex structure carried by $P$. Then each embedded surface carried by $P$ has minimal genus in its homology class.
\end{teo}

{\bf Acknowledgements.}
We wish to express our gratitude to Stephane Baseilhac, Riccardo Benedetti, Paolo Lisca, Dylan Thurston and Vladimir Turaev for their interest and criticism on this work.
\section{Branched shadows of $4$-manifolds}

In this section we recall the notions of shadow of $4$-manifold, of branching and the main results on this topic.
For a complete account on shadows see \cite{Tu} and \cite{Tu2}; for an introductory one, see \cite{Co5}; for a detailed account of branched shadows we refer to
\cite{Cobr}. From now on, all the manifolds and homeomorphisms will be smooth unless explicitly stated. 

\subsection{Shadows of $4$-manifolds}\label{basicfacts}
\begin{figure}[h!]
  \centerline{\includegraphics[width=8.4cm]{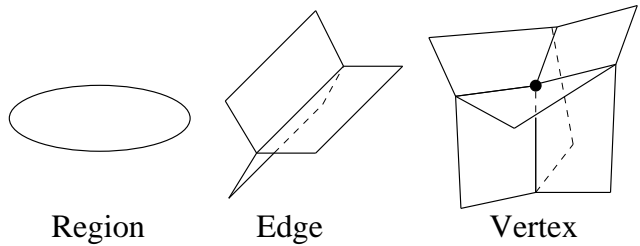}}
  \caption{The three local models of a simple polyhedron. }
  \label{fig:singularityinspine}
\end{figure}

\begin{defi} A simple polyhedron $P$ is a finite $2$-dimensional CW complex whose
local models are those depicted in Figure
\ref{fig:singularityinspine}; the set of points whose neighborhoods have models
of the two rightmost types is a $4$-valent graph, called
{\it singular set} of the polyhedron and denoted by $Sing(P)$. The
connected components of $P-Sing(P)$ are the {\it regions} of $P$.
A simple polyhedron whose singular set is connected
and whose regions are all discs is called {\it
standard}. 
\end{defi}
\begin{defi}[Shadow of a $4$-manifold]\label{def:flat}
  Let $M$ be a smooth, compact and oriented $4$-manifold.  $P \subset M$ is a
  \emph{shadow} for $M$ if:
  \begin{enumerate}
  \item $P$ is a closed polyhedron  embedded in $M$ so that $M-P$ is diffeomorphic to $\partial M\times (0,1]$;
  \item $P$ is \emph{flat} in $M$, that is for each point
  $p\in P$ there exists a local chart $(U,\phi)$ of $M$ around $p$ such
  that $\phi(P\cap U)$ is contained in $\mr^3 \subset \mr^4$ and in this chart the
pair $(\mathbb{R}^3\cap \phi (U),\mathbb{R}^3\cap \phi(U\cap P))$ is diffeomorphic
to one of the models depicted in Figure \ref{fig:singularityinspine}.
\end{enumerate}
\end{defi}

From now on, all the embedded polyhedra will be flat unless explicitly stated.
\begin{rem}
Note that the original definition of shadows was given in the PL setting by Turaev in \cite{Tu}
. But in four dimensions the smooth and the PL setting are equivalent, that is for each PL-structure on a compact manifold there exists a unique compatible (in a suitable sense) smooth structure. Definition \ref{def:flat} is the natural translation of the notion of shadow to the smooth setting. 
\end{rem}

A necessary and sufficient
condition (see \cite{Co}) for a $4$-manifold $M$ to admit a shadow is that $M$ is a $4$-handlebody, that is $M$ admits a handle decomposition without $3$ and $4$-handles.
In particular, $\partial M$ is a non empty connected $3$-manifold.
From now on, all the manifolds will be $4$-handlebodies unless explicitly stated. 

\begin{figure}
  \centerline{\includegraphics[width=8.4cm]{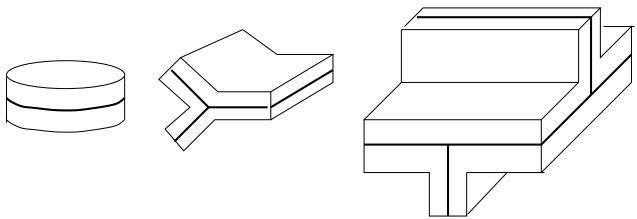}}
  \caption{The three type of blocks used to thicken a polyhedron to a $3$-manifold.  }
  \label{fig:spineblocks}
\end{figure}

Suppose that $P$ is a surface embedded in a oriented
$4$-manifold $M$. In general we cannot reconstruct the tubular
neighborhoods of $P$ by using only its topology, since their structure depends on the self-intersection number of $P$ in $M$.  To encode the topology of a
neighborhood of $P$ in $M$ we need to ``decorate" $P$ with
some additional information; when $P$ is an oriented surface,
the Euler number of its normal bundle is a sufficient datum.

We describe now the basic decorations we need for a general polyhedron $P$.
Let us denote $\mathbb{Z}[\frac{1}{2}]$ the group of integer multiples of $\frac{1}{2}$. There are two canonical {\it colorings} on the regions of $P$, i.e. assignments of elements of $\mathbb{Z}_2$ or $\mathbb{Z}[\frac{1}{2}]$, the second depending on a flat embedding of $P$ in an oriented  $4$-manifold. They are:

{\bf The $\mathbb{Z}_2$-gleam of $P$}, constructed as follows.  Let $D$ be the (open) $2$-cell associated to a given region of $P$ and $\overline{D}$ be
the natural compactification $\overline{D}=D\cup S^1$ of the (open) surface represented by
$D$ and $cl(D)$ the closure of $D$ in $P$. The embedding of $D$ in $P$ extends to a map
$i:\overline{D}\to cl(D)$ which is injective in $int(\overline{D})$,
locally injective on $\partial \overline{D}$ and which sends
$\partial \overline{D}$ into $Sing(P)$. Using $i$ we can
``pull back" a small open neighborhood of $D$ in $P$ and
construct a simple polyhedron $N(D)$ collapsing on $\overline{D}$
and such that the map $i$ extends as a local homeomorphism
$i':N(D)\to P$ whose image is contained in a small
neighborhood of the closure of $D$ in $P$. When $i$ is an embedding of $\overline{D}$ in $P$, then $N(D)$
turns out to be homeomorphic to a neighborhood of $D$ in $P$
and $i'$ is its embedding in $P$. In general, $N(D)$ has the following structure: each boundary
component of $\overline{D}$ is glued to the core of an annulus or
of a M\"obius strip and some small discs are glued along half of
their boundary on segments which are properly embedded in these
annuli or strips and cut transversally once their cores. We
define the $\mathbb{Z}_2$-gleam of $D$ in $P$ as the
reduction modulo $2$ of the number of M\"obius strips used to
construct $N(D)$. This coloring only depends on
the combinatorial structure of $P$.

\begin{figure}
  \centerline{\includegraphics[width=6.4cm]{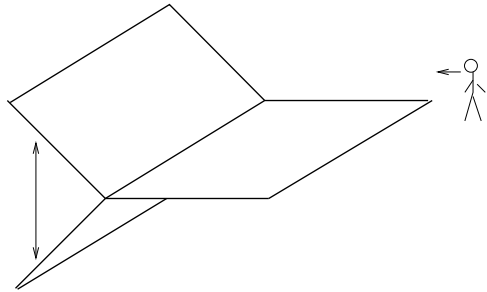}}
  \caption{The picture sketches the position of the polyhedron in a
    $3$-dimensional slice of the ambient $4$-manifold.  The direction
    indicated by the vertical double arrow is the one along which the
    two regions touching the horizontal one get separated.  }
  \label{fig:divergingdirection}
\end{figure}

{\bf The gleam of $P$}\label{embpoly}, constructed as follows. Let us now suppose that $P$ is flat in an oriented $4$-manifold $M$, with $D$, $\overline{D}$, $cl(D)$ and $i:\overline{D}\rightarrow P$ as above. 
Orient $D$ arbitrarily and orient its normal bundle $n$ in $M$ so that the orientation on the global space of $n$ coincides with that of $M$.
Pulling back $n$ to $\overline{D}$ through $i$, we get an oriented disc bundle over $\overline{D}$ we will call the ``normal bundle" of $\overline{D}$; we claim that the projectivization of this bundle comes equipped with a section defined on $\partial \overline{D}$. 
Indeed, since $P$ is locally flat in $M$, for each point $p\in \partial \overline{D}$, if $i(p)\in Sing(P)$ there exists a smooth $3$-ball $B_{i(p)}\subset M$ around $i(p)$ in which $P$ appears as in Figure \ref{fig:divergingdirection}. Then, the intersection in $T_{i(p)} M$ of the ($3$-dimensional) tangent bundle of $B_{i(p)}$ with the normal bundle to $cl(D)$ in $i(p)$ gives a normal direction to $cl(D)$ in $i(p)$ (indicated in the figure) whose pull-back through $i$ is the seeked section in $p$. Hence a section of the projectivized normal bundle of $\overline{D}$ is defined on all $\partial \overline{D}$: we then define $gl(D)$ be equal to $\frac{1}{2}$ times the obstruction to extend this section to the whole $\overline{D}$; such an obstruction is an element of $H^2(\overline{D},\partial \overline{D};\pi_1(S^1))$, which is canonically identified with $\mathbb{Z}$ since $M$ is oriented.
Note that the gleam of a region is integer if and only if its $\mathbb{Z}_2$-gleam is zero. 

Using the fact
that the $\mathbb{Z}_2$-gleam is always defined, Turaev 
generalized~\cite{Tu} the
notion of gleam to non-embedded polyhedra as follows:
\begin{defi}
A {\it gleam} on a simple polyhedron $P$ is a coloring on the regions of
$P$ with values in $\mathbb{Z}[\frac{1}{2}]$ such that the color of
a region is integer if and only if its $\mathbb{Z}_2$-gleam is
zero.
\end{defi}
\begin{teo}[Reconstruction Theorem \cite{Tu}]\label{teo:reconstruction}
Let $P$ be a polyhedron with gleams~$gl$; there exists a
canonical reconstruction map associating to $(P,gl)$ a pair $(M_P,P)$
where $M_P$ is a smooth, compact and oriented $4$-manifold, and $P\subset M$ is a shadow of $M$ (see Definition \ref{def:flat}).
If $P$ is a polyhedron flat in a smooth and oriented  $4$-manifold and $gl$ is the gleam of $P$ induced by its embedding, then $M_P$ is diffeomorphic to a compact neighborhood of $P$ in $M$.
\end{teo}
The proof is
based on a block by block reconstruction procedure similar to the
one used to describe $3$-manifolds by means of their spines.
For instance, if $P$ is a standard polyhedron, for each of the three local patterns of Figure \ref{fig:singularityinspine}, we consider the $4$-dimensional thickening given by the
product of an interval with the corresponding $3$-dimensional
block shown in Figure \ref{fig:spineblocks}. All
these thickenings are glued to each other according to the
combinatorics of $P$ and to its gleam.

By Theorem \ref{teo:reconstruction}, to study $4$-manifolds one can either
use abstract polyhedra equipped with
gleams or embedded polyhedra. The latter approach is more
abstract, while the former is purely combinatorial; 
we will use both approaches in
the following sections. The translation in the combinatorial setting of the definition of shadow of a $4$-manifold is the following:

\begin{defi}[Combinatorial shadow]
A polyhedron equipped with gleams $(P,gl)$ is said to be a {\it
shadow} of the $4$-manifold $M$ if $M$ is diffeomorphic to the manifold associated to $(P,gl)$ by means of
the reconstruction map of Theorem
\ref{teo:reconstruction}.
\end{defi}

\subsection{Branched shadows}

From now on, for the sake of simplicity all the
polyhedra will be standard without
explicit mentioning. Given a polyhedron $P$ we define the notion of
{\it branching} on it as follows:
\begin{defi}[Branching condition]\label{branchingcondition}
A branching $b$ on $P$ is a choice of an orientation for each
region of $P$ such that for each edge of $P$, the
orientations induced on the edge by the regions containing it do
not coincide. \end{defi}
\begin{rem}
This definition corresponds to the definition of ``orientable
branching'' of \cite{BP}.
\end{rem}

We say that a polyhedron is {\it branchable} if it admits a
branching and we call {\it branched polyhedron} a pair
$(P,b)$ where $b$ is a branching on $P$.
Each edge $e$ of $Sing(P)$ is touched locally by three regions (which at a global level could be be non-distinct). The branching condition implies that these three oriented regions will induce two-times one orientation on $e$ and once the opposite one. We say that $e$ is oriented {\it according to the branching} of $P$ if it is equipped with the former orientation.

\begin{defi}
Let $(P,gl)$ be a shadow of a $4$-manifold $M$; $P$ is
said to be {\it branchable} if the underlying polyhedron is. We
call {\it branched shadow} of $M$ the triple $(P,gl,b)$ where
$(P,gl)$ is a shadow and $b$ is a branching on $P$. 
When this will not cause any confusion, we will not
specify the branching $b$ and we will simply write $(P,gl)$.
\end{defi}
\begin{prop}[\cite{Cobr}]\label{prop:existsbranchedshadow}
Any $4$-manifold admitting a shadow admits also a branched shadow.
\end{prop}

\begin{figure} [h!]
   \centerline{\includegraphics[width=11.4cm]{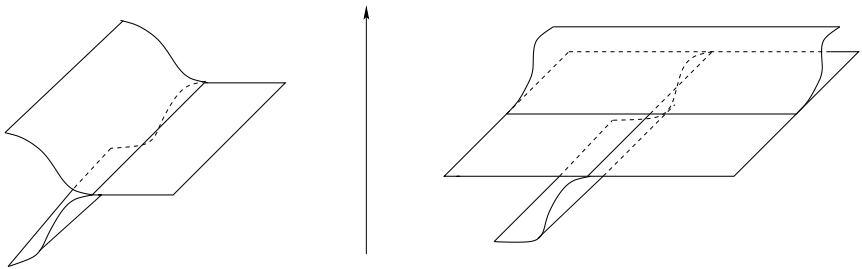}}
   \caption{How a branching allows a smoothing of the polyhedron: the
  regions are oriented so that their projection on the ``horizontal" (orthogonal to the drawn vertical direction) plane is orientation preserving.}\label{branching}
\end{figure}

A branching on a shadow allows us to smoothen its singularities
and equip it with a smooth structure as shown in Figure
\ref{branching}. This smoothing can be performed also inside the
ambient manifold obtained by thickening the shadow; the shadow locally appears as in 
Figure \ref{branching}, where the two regions orienting the edge
in the same direction approach each other so that, for any auxiliary riemannian metric on the ambient manifold, all the derivatives of their
distance go to zero while approaching the edge.
\begin{defi}\label{defi:preferred}
Let $e$ be an edge of a branched polyhedron $P$ and let $R_i$,$R_j$ and $R_k$ be the regions of $P$ containing it in their boundary. Then $R_i$ is said to be the {\it preferred region of $e$} if it induces the opposite orientation on $e$ with respect to those induced by $R_j$ and $R_k$.
\end{defi}

\begin{figure}
   \centerline{\includegraphics[width=6.4cm]{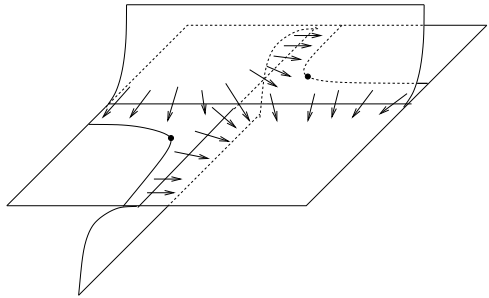}}

  \caption{ In this figure we show how the maw behaves near a
  vertex of a branched polyhedron.
  }\label{fig:maw}
\end{figure}

The following proposition is a consequence of the fact that a $4$-handlebody retracts on its shadows.
\begin{prop} \label{prop:branched presentation}
Let $P$ be a branched shadow of a $4$-manifold
$M$, and let $R_i,\ i=1,\ldots,n$ and $e_j,j=1,\ldots,m$ be
respectively the regions and the edges of $P$ oriented according
to the branching of $P$. Then $H_2(M;\mz)$ is the kernel of
the boundary application $\partial :\mz^n \to
\mz^m$. Moreover
$H^{2}(M;\mz)$ is the abelian group generated by
the cochains $\hat{R}_i,\ i=1,\ldots,n$ dual to the regions of $P$
subject to the relations generated by the coboundaries of the edges having the form  $\delta(e_j)=-\hat{R}_i+\hat{R}_j+\hat{R}_k$ where $R_i$ is the preferred region of $e_j$. Analogously, $H^1(M;\mathbb{Z})$ is the abelian group generated by the cochains $\hat{e}_1,\dots, \hat{e}_m$ dual to the edges subject to the relations generated by the coboundaries of the vertices having the form $\delta(\hat{v}_s)=\hat{e}_t+\hat{e}_u-\hat{e}_v-\hat{e}_z$ where $e_t$ and $e_u$ (resp. $e_v$ and $e_z$) are the edges whose final (initial) endpoint is $v_s$.
\end{prop}
Given a shadow $(P,gl)$ of a $4$-manifold $M$, there are three cochains representing classes in $H^2(M;\mathbb{Z})$ naturally associated to $(P,gl)$. 
The first one is the {\it Euler cochain} of $P$, denoted $Eul(P)$ and constructed as follows.
Let $m$ be the vector field tangent to $P$ (using the smoothed structure given by the branching)
which near the center of the edges points inside the preferred
regions; we extend $m$ in a neighborhood of the
vertices as shown in Figure \ref{fig:maw}.

The field constructed above is the {\it maw}
of $P$. For each region $R_i$ of $P$, the maw
gives a vector field defined near $\partial R_i$, so
it is possible to extend this field to a tangent field on
the whole $R_i$ having isolated singularities of indices
$\pm 1$; let $n_i$ be the algebraic sum of these indices over the
region $R_i$. The Euler cochain is defined as $Eul(P)=\sum_i
n_i\hat{R}_i$; its meaning will be clarified in the next subsection.

\subsection{Branched shadows and almost complex structures}\label{sub:brshalmostcomplex}

As before, let $M$ be an oriented $4$-handlebody and let $g$ be a fixed auxiliary riemannian metric on $M$. In this subsection we recall how a
branched shadow determines a pair $(M,[J])$ where
$[J]$ is a homotopy class of almost complex structures on $M$ suitably compatible with $g$,
and recall that for each such class, there exists a branched shadow of
$M$ encoding it. 

\begin{defi}
An almost complex structure $J$ on an oriented $4$-manifold $M$ is a smooth morphism $J:TM\to TM$ such that for each point $p\in M$ it holds $J^2=-Id$. We say that $J$ is {\it positive} if at each point $p\in M$ there exists a positive (with respect to the orientation of $M$) basis of $T_pM$ of the form $(x,J(x),y,J(y))$. We say that $J$ is orthogonal with respect to $g$ if for each $p\in M$ the map $J:T_pM\to T_pM$ is $g$-orthogonal.
\end{defi}

Let $P$ be a branched shadow of $M$. In each of the local blocks used to reconstruct
$M$ from $P$ as in Theorem \ref{teo:reconstruction} (these blocks are the products of an interval with the
$3$-dimensional blocks of Figure \ref{fig:spineblocks}), $P$ is smoothly embedded in a non symmetrical way (see Figure \ref{branching}). As in Figure \ref{branching}, we 
choose an
{\it horizontal} $2$-plane and
orient it according to the branching of $P$, so that each of
the basic building blocks of $M$ is equipped with a
distribution of oriented $2$-planes, denoted as $T(P)$.

Let $V(P)$ be the field of
oriented {\it vertical} $2$-planes of $TM$ which are pointwise positively 
$g$-orthogonal to the planes of $T(P)$.
We define an almost complex structure $J_P$ by requiring that its restriction to the two fields $T(P)$ and $V(P)$ acts pointwise as a $\frac{\pi}{2}$ positive rotation and we extend this action linearly to the whole $TM$. By construction $J_P$ is positive and $g$-orthogonal. 
\begin{teo}[\cite{Cobr}]\label{teo:surjectivityofrefinedreconstruction}
The restriction to branched shadows of the reconstruction map of Theorem \ref{teo:reconstruction} is a well defined surjection on the pairs $(M,[J])$ with $[J]$ homotopy class of positive almost complex structures on $M$.  
\end{teo}
Our construction splits the tangent bundle of $M$ as the sum of two
linear complex bundles $V(P)$ and $T(P)$, hence the first Chern
class of $TM$, viewed as a complex bundle using the almost complex
structure $J_P$, is equal to
$c_1(T(P))+c_1(V(P))$. The following proposition
(whose proof is identical to that of Proposition 7.1.1 of \cite{BP}) is a recipe to calculate $c_1(T(P))$:
\begin{prop}\label{prop:eulercochain2}
The class in $H^2(P;\mz)$ represented by the Euler cochain $Eul(P)$
coincides with the first Chern class of the horizontal plane field
$T(P)$ of $P$ in $M$.
\end{prop}
To calculate $c_1(V(P))$, we first need to fix some notation.
Note that a branching on a polyhedron allows one to define a set of $1$-cochains with values in $\mathbb{Z}_2$ we will call ``Up\&Down"-cochains as follows. In a neighborhood of each vertex, we can canonically identify $2$ regions out of the $6$ surrounding it: those which are not horizontal in Figure \ref{branching}. Let us choose arbitrarily around each vertex of $P$ which of these two regions is ``up" and which is ``down"; each edge touching a vertex is contained in one of these regions. We then define a $\mathbb{Z}_2$-cochain on an edge to be equal to $1$ iff, following the up (down) region containing the edge in one of its endpoints, we get to the down (up) region in the other endpoint. Our definition of the values of the cochain depends of course on the choices of the up and down regions around the vertices; one can check that the difference of two ``Up\&Down"-cochains is a coboundary (see Proposition \ref{prop:branched presentation}).
\begin{lemma}
Let $P$ be a branched polyhedron and $ud$ be an ``{\rm Up\&Down}" cochain on $P$. Then $\delta ud=gl_2$ that is, the coboundary of $ud$ is the $\mathbb{Z}_2$-cochain representing the $\mathbb{Z}_2$-gleam. In particular the $0$ cochain is an ``{\rm Up\&Down}"-cochain iff $P$ is a spine of a $3$-manifold.  
\end{lemma}
\begin{defi}
A {\it $\mathbb{Z}[\frac{1}{2}]$-lift} of an ``Up{\rm \&}Down"-cochain $ud$ is a $1$-cochain $UD$ with values in $\mathbb{Z}[\frac{1}{2}]$ such that the evaluation of $UD$ on an edge $e$ is an odd multiple of $\frac{1}{2}$ iff $ud(e)=1$. The {\it canonical lift} of an ``Up{\rm \&}Down"-cochain $ud$ is the $\mathbb{Z}[\frac{1}{2}]$-lift whose value on an edge $e$ is $\frac{1}{2}$ if $ud(e)=1$ and is $0$ otherwise.
\end{defi}

We are now ready to define the {\it gleam cochain} $gl(P)$ as
$gl(P)=\sum_i gl(R_i)\hat{R}_i$, where $\hat{R}_i$ is the cochain
dual to $R_i$, the coefficient $gl(R_i)$ is the gleam of the
region $R_i$ and the sum ranges over all the regions of $P$.
Note that since the gleams can be half-integers, it is not a priori obvious that the gleam cochain represents an integer class in $H^2(M;\mathbb{Z})\cong H^2(P;\mathbb{Z})$; to by pass this problem, 
it can be checked that $gl(P)=\frac{1}{2}\delta(UD)+gl_{\mathbb{Z}}$, where $\delta(UD)$ is the coboundary of the canonical lift of an Up\&Down cochain and $gl_{\mathbb{Z}}$ is an integer coefficient cochain. Since all the possible choices of the $1$-cochain $UD$ differ by a coboundary, the coefficients of both the cochain $\delta(UD)$ and the cochain $gl_{Z}$ are well defined. Hence $gl(P)([S])=gl_{\mz}([S])$ for each cycle $[S]\in H_2(P;\mz)$ and often, when considering cohomology classes, we will refer to $gl(P)$ instead of $gl_{\mz}(P)$: \begin{lemma}[\cite{Cobr}]\label{glclass}
The gleam cochain $gl(P)$ represents in $H^2(M;\mathbb{Z})$ the
first Chern class of the field of oriented $2$-planes $V(P)$.
\end{lemma}
\begin{cor}\label{Chern}
The first Chern class of the almost complex structure on $M$ associated
with $P$ is represented by the cochain $c_1(P)=Eul(P)+gl(P)$.
\end{cor}

\subsection{Branched shadows and complex structures}\label{brshcomplex}
\begin{defi}
An almost complex structure $J$ on a smooth $4$-manifold $M$ is said
to be {\it integrable} or {\it complex} if for each point $p$ of $M$
there is a local chart of $M$ with values in $\mathbb{C}^2$ transforming $J$ into the complex structure of $\mathbb{C}^2$.
\end{defi}

Let $P$ be a branched shadow embedded in a complex $4$-manifold $M$. Up to perturbing the embedding of $P$ through a small isotopy we can suppose that there is only a finite
number of points $p_1,\ldots, p_n$ and $q_1,\ldots,q_m$, contained in the regions of $P$ where $T_{p_i}P$ (resp. $T_{q_j}P$) is a complex plane such that the orientations induced by the branching of $P$ and by the complex structure of $M$ coincide (resp. do not coincide).

\begin{defi}
The points $p_1,\ldots,p_n$ are called
{\it positive complex} points of $P$ or simply {\it positive}
points. Analogously, the points $q_1,\ldots,q_m$ are called {\it
negative complex} points of $P$ or {\it negative} points. All the
other points of $P$ are called {\it totally real}.
\end{defi}

To each complex point $p$ of a region $R_i$ of $P$ we can assign an integer number
called its {\it index}, denoted $i(p)$, as follows. Fix a small disc $D$ in $R_i$ containing $p$ and no
other complex point and let $N$ be the radial vector field
around $p$. The field $J(N)$ on
$\partial D$ is a vector field transverse to $P$ since no point on
$D-p$ is complex. Let $\pi (J(N))$ be the projection of this field
onto the normal bundle of $D$ in $M$. Since $D$ is contractile,
this bundle is trivial and we can count the number $\nu(p)$
of twists performed by $\pi(J(N))$ while following $\partial D$
($D$ and $M$ are oriented). The index of $p$ is: $i(p)=\nu(p)+1$. 
Moreover, we define $\nu(R_i)$ as the sum
over all the complex points $p$ of $R_i$ of $\nu(p)$. 

Up to a small perturbation by an isotopy of the embedding of $P$ in $M$ we can assume that all the indices of the complex points of $P$ are equal to $\pm 1$.

\begin{defi}
A complex point $p$ of $P$ whose index is equal to $1$ is {\it elliptic}, if its index is $-1$ it is {\it hyperbolic}.
\end{defi}

We define the index $c_1(p)$ associated to each complex point $p$ of a region $R_i$ as follows. Let $D$ and $N$ be as above; complete $N$ on $\partial D$ to a basis of $TD$ by using
the field $T=T\partial D$ tangent to the boundary of $D$. The pair
of fields $(N,T)$ gives a basis of $TD$ in each point $q$
of $\partial D$, moreover, since no such point is complex, they can
be completed to a positive complex basis of $T_qM$ given by 
$(N,J(N),T,J(T))$. Let now $\frac{\partial}{\partial z}$ and
$\frac{\partial}{\partial w}$ be two vector fields defined on a
neighborhood of $D$ in $M$ such that $(\frac{\partial}{\partial
z},\frac{\partial}{\partial w})$ is pointwise a complex basis of
$TM$. Then, on each point $q$ of $\partial D$ we can compare the
two complex bases given by $(N+J(N),T+J(T))$ and
$(\frac{\partial}{\partial z},\frac{\partial}{\partial w})$ by
considering the determinant $det_q$ of the change of basis from
the latter to the former basis. The value of the index of $det_q$ around $0$ in
$\mathbb{C}$ while $q$ runs across $\partial D$ according to the
orientation of $D$, is defined to be $c_1(p)$. We define $c_1(R_i)$ as the sum of $c_1(p)$ over all the complex points $p$ of $R_i$.

The following is a straightforward generalization of a result proved by Bishop \cite{Bi}, Chern and Spanier \cite{CS} and Lai \cite{La}, in the case of closed surfaces: 
\begin{teo} \label{teo:bishop}
Let $R_i$ be a surface with boundary contained in complex manifold $M$ such that $\partial R_i$ does not contain complex points and let $I^+=\sum_{i=1,\ldots,n}i(p_i)$ and
$I^-=\sum_{j=1,\ldots,m}i(q_j)$. Then the following equalities
hold:
$$I^+=\frac{1}{2}(\chi (R_i)+\nu(R_i)+c_1(R_i))$$
$$I^-=\frac{1}{2}(\chi(R_i)+\nu(R_i)-c_1(R_i))$$
\end{teo}

The following is the analogue in the world of shadows of the well known result of V.M. Harlamov and Y. Eliashberg on annihilation of pair of complex points in real surfaces embedded in complex manifolds:
\begin{lemma}[\cite{Cobr}]\label{lem:branchedHE}
Let $R_i$, $R_j$ and $R_k$ be three regions of $P$ adjacent along
a common edge $e\in Sing(P)$ so that $R_i$ is the preferred region of $e$. 
There exists an isotopy
$\phi_t:P\to M,\ t\in [0,1]$ whose support is contained in a small
ball $B$ around the center of $e$ such that $\phi_1(B\cap P)$ contains three
more complex points $p_i$, $p_j$ and $p_k$ respectively in $R_i$,
$R_j$ and $R_k$ whose indices are respectively $\nu(p_i)=\pm 1$, $\nu(P_j)=\mp 1
$ and $\nu(P_k)=\mp 1$.
\end{lemma}

The above lemma suggests the following:
\begin{defi}\label{def:ipluseiminus}
The {\it positive index} and {\it negative index} cochains of $P$, denoted respectively $I^+(P)$ and $I^-(P)$ are the $2$-cochains given by $\Sigma_i I^{\pm}(R_i)\hat{R}_i$, where $i$ ranges over all the regions of $P$.
\end{defi}
\begin{teo}[\cite{Cobr}]
The cohomology classes $[I^{\pm}(P)]\in H^2(M;\mathbb{Z})$ are invariants of the embedding of $P$ in $M$ up to isotopy.
\end{teo}
\begin{prop}[\cite{Cobr}]\label{prop:imenougualediff}
The almost complex structure carried by $P$ is homotopic to the ambient complex structure iff $I^-(P)=0 \in H^2(M;\mathbb{Z})$.
\end{prop}

\section{Branched shadows and Stein domains}
\subsection{Classical facts on Stein domains}
We summarize here the main facts on Stein domains, see Chapter XI of \cite{GS} for a detailed account.
\begin{defi}[Stein domain]\label{defi:stein}
A {\it Stein domain} is a compact,complex $4$-manifold $M$ equipped with a smooth Morse function $\phi:M\rightarrow [0,1]$ such that $\partial M=\phi^{-1}(1)$ and $\phi$ is strictly plurysubharmonic, that is, its complex Hessian is a positive definite Hermitian form.
\end{defi}
A $4$-manifold diffeomorphic to $D^4\cup 1$-handles admits a canonical Stein domain structure called the {\it standard structure} (see \cite{El}).
With the notation above, each preimage $\phi^{-1}(t)$ of a regular value of $\phi$ is equipped with a {\it contact structure} $\alpha_t$ that is a totally non-integrable distribution of oriented $2$-planes. These $2$-planes are the only complex tangent planes of $TM$ entirely contained in $T(\phi^{-1}(t))$.

\begin{defi}[Legendrian link]\label{defi:legendrian}
A link $L$ in an oriented $3$-manifold $N$ equipped with a contact structure is {\it Legendrian} if it is everywhere tangent to the contact structure. The vector field in $N$ which is positively transverse to the contact structure equips each Legendrian link with a framing called {\it Thurston-Bennequin framing}, and denoted $tb$.  
\end{defi}
Each link in a $3$-manifold equipped with a contact structure is isotopic through a $C^0$-small isotopy to a Legendrian one. The following is the well known theorem of Y. Eliashberg and R. Gompf (see \cite{El} and \cite{Go}).
\begin{teo}[Legendrian surgery]\label{teo:gompf-eliashberg}
A smooth, compact, connected, oriented $4$-manifold $M$ admits a Stein domain structure (inducing the given orientation) if and only if it has a handle decomposition such that the $2$-handles are attached to the union of $0$ and $1$-handles equipped with its standard structure along Legendrian knots and the framing coefficient on each knot is $tb-1$. 
\end{teo}
Through a small isotopy, one can always move a Legendrian knot and decrease its $tb$ framing by a unit. This can be performed in two ways, which in the standard way of drawing Legendrian knots in $\mr^3$ (see \cite{GS}) correspond to adding a positive or negative ``zig-zag". This implies that in the above theorem one can require the framing to be $tb+k$, $k\leq -1$. 
\begin{rem}\label{rem:hyperbpoints}
The construction underlying Theorem \ref{teo:gompf-eliashberg} provides a handle decomposition such that the core of the $2$-handle is totally real w.r.t. the ambient complex structure. When one applies the theorem with framing $tb+k$ with $k<-1$, the core of the handle (oriented arbitrarily) contains $-1-k$ hyperbolic complex points and no elliptic complex points. One can choose to have a $h^+\geq 0$ hyperbolic positive points and $h^-\geq 0$ hyperbolic negative points provided that $h^++h^-=-1-k$. This corresponds to choosing how many positive and negative ``zig-zags" add to the Legendrian knot before applying Theorem \ref{teo:gompf-eliashberg}. In particular, when $k<-1$ there are in principle more than one Stein domains structures, since their almost complex structures have different $c_1$ (see Theorem \ref{teo:bishop}) on the core of the handle. These structures could anyway be isotopic due to some global property of the ambient manifold. 
\end{rem}

\subsection{Branched shadows of Legendrian curves} \label{legendriancurves}

We now recall Polyak's method for drawing Legendrian links in the boundary of a Stein domain collapsing over a surface (for further details see \cite{Po}), and then construct a branched polyhedron associated to each such link embedded in the Stein domain.

Let $\Sigma$ be an oriented surface and let $M_{\Sigma}$ be the total space of the disc bundle over $\Sigma$ with Euler number equal to $-\chi(\Sigma)$. It is a standard fact that $M_{\Sigma}$ can be equipped with a Stein domain structure with respect to which $\Sigma$ is totally real and whose Morse function $\phi:M_{\Sigma}\to [0,1]$ is equal to the distance from $\Sigma$ (in a fixed riemannian metric on $M_{\Sigma}$).  The above facts hold also when $\partial \Sigma\neq 0$ with two differences: $M_{\Sigma}$ is a neighborhood of $\Sigma$ in a Stein domain in which $\Sigma$ is totally real and $M_{\Sigma}$ is cut by a function $\phi$ which coincides with the distance from $\Sigma$ near the interior part of $\Sigma$ and is slightly perturbed near its boundary (see \cite{Fo} and \cite{NW} for further details). 

As showed in \cite{Po}, projecting into $\Sigma$ a Legendrian knot $k$ contained in $\partial M_{\Sigma}$  produces a ``wavefront", that is a smooth curve $c$ possibly containing some cusps and equipped with a {\it coorientation} that is, a normal vector field $v$. Moreover, one can recover the knot from the decorated curve. Indeed, in a local chart $(x_1,x_2,x_3,x_4)$ around $\Sigma=(x_1,x_2)$ where the complex structure of $M_{\Sigma}$, looks like $J(\frac{\partial}{\partial x_1})=\frac{\partial}{\partial x_3}$ and $J(\frac{\partial}{\partial x_2})=-\frac{\partial}{\partial x_4}$ the knot $k$ is given by $c+J(v)$, where, for instance, by $\vec p+\alpha \frac{\partial}{\partial x_1}, \alpha \in \mr$ we denote the point whose coordinates are $\vec p+(\alpha,0,0,0)$.

We now associate to each oriented Legendrian link $l$ in $\partial  M_{\Sigma}$ a branched polyhedron embedded in $M_{\Sigma}$ whose boundary contains $l$ as follows. 
Suppose first that the projection of $l$ in $\Sigma$ does not contain cusps; then the mapping cylinder $P_l$ of the projection of $l$ in $\Sigma$ will do: orient the regions coming from $\Sigma$ as subsets of $\Sigma$ and the (annular) regions coming from the components of $l$ so that they induce on the projection of $l$ the positive orientation. This defines a branching on $P_l$ and allows us to smoothen its singularities in $M_{\Sigma}$ as explained in Figure \ref{branching}. It can be checked that $P_l$ is totally real and that, while smoothing the singularities, the singular set (which initially coincides with the projection of $l$) is slightly pushed at its right (both $\Sigma$ and $l$ are oriented).

\begin{figure}
   \centerline{\includegraphics[width=10.4cm]{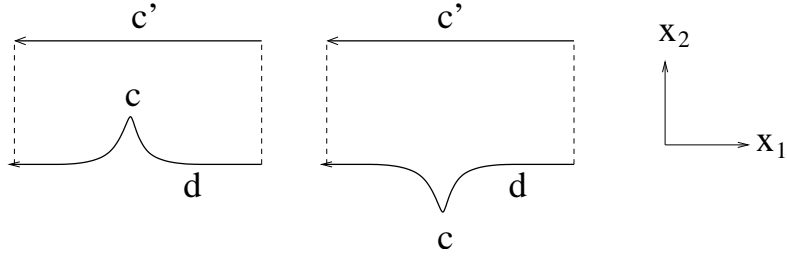}}
  \caption{Smoothing a cusp.
  }\label{fig:smoothcusp}
\end{figure}

Let us now study how to construct $P_l$ in a neighborhood of a cusp; we say that a cusp $c$ is {\it positive} if it points in the left of $c$ and {\it negative} otherwise: indeed recall that both the surfaces and the curves we are dealing with are oriented. For instance, the leftmost drawing in Figure \ref{fig:smoothcusp} is a negative cusp and the central one is a positive cusp. Let us fix a local chart around the cusp in which $\Sigma =(x_1,x_2)$, let $d$ be the arc of $l$ whose projection forms $c$ and $c'$ be a straight oriented arc in $\Sigma$ running parallel to the cusp and lying at its right (see Figure \ref{fig:smoothcusp}).

\begin{lemma}
With the notation above, there exists a smooth disc $D(c)$ in the local chart around $c$ such that $D(c)\cap \Sigma=\partial D(c)\cap \Sigma=c'$, $D(c)\cap d=\partial D(c)\cap d=d$ and such that $\partial D(c) - d\cup c'$ is formed by two dotted arcs of Figure \ref{fig:smoothcusp}; moreover the disc is tangent to $\Sigma$ along $c'$ and it induces the positive orientation on $c'$. If $c$ is positive then $D$ is totally real regardless of the coorientation of $c$, if $c$ is negative then $D(c)$ contains a hyperbolic complex point (either positive or negative according to the coorientation of $c$).
\end{lemma}\label{lem:cuspanalisis}
\begin{rem}
It is worth to note that the asymmetry in the above statement is due to the fact that $D(c)$ is chosen to be tangent along $c'$ to $\Sigma$ and so, changing the orientation of $d$ and hence the sign of the cusp, changes also the position of $D(c)$, not only its orientation.  
\end{rem}
 \begin{prf}{2}{
Fix a local chart in which $\Sigma=(x_1,x_2)$, $J(\frac{\partial}{\partial x_1})=\frac{\partial}{\partial x_3}$ and $J(\frac{\partial}{\partial x_2})=-\frac{\partial}{\partial x_4}$.
We first exhibit $D(c)$ and then show how to count the number of its complex points.
Fix an orthogonal vector field $v=(v_1(s),v_2(s),0,0)$ on $c(s),\ s\in [-1,1]$ and consider the rectangle $R$ tangent to $\mr^2$ along $c'$, parametrized as follows (using the coordinates of Figure \ref{fig:smoothcusp}): $R(s,t)=(s,t,0,0)+\frac{1}{2}t^2(J(v(s))),\ (s,t)\in [-1,1]\times [-1,0]$. The boundary curve $R(s,-1)$ is parallel in $\mr^4-\Sigma$ to $d$; let $R'$ be the rectangle forming the parallelism. We can choose $R'$ so that the disc $D(c)=R\cup R'$ is smooth and totally real at its boundary. 
We count the number of complex points of $D(c)$ using Theorem \ref{teo:bishop}. To do this, we need non zero sections respectively of the normal bundle of $D(c)$ in $\mr^4$ and of the determinant line bundle of $(\mr^4,J)$ restricted to $D(c)$. It can be checked that the field $B$ along $\partial D(c)$ defined as follows extends to a non zero section of $D(c)$. On $c'$ let $B=J(v(s))$, on $d$ let $B=v(s)$ and on the two dotted edges of Figure \ref{fig:smoothcusp}, let $B=cos(\theta) J(v(s))+sin(\theta) v(s)$ where $\theta \in [0,\frac{\pi}{2}]$ is a parameter for the edges oriented according to the orientation induced from $D(c)$ on them. More easily, a non zero section of the determinant line bundle of $(\mr^4,J)$ is given by $\frac{\partial}{\partial x_1}\wedge \frac{\partial}{\partial x_2}$.
Now we need to calculate the $\nu$ and $c_1$ as explained in Subsection \ref{brshcomplex}; the image through $J$ of the field $N$ tangent to $D(c)$ and pointing outside it along its boundary has the following behavior: on $c'$ it is equal to $-\frac{\partial}{\partial x_4}$, on $d$ it is $-v(s)$ and on the vertical edges it is $\pm \frac{\partial}{\partial x_3}$. We now follow $\partial D(c)$ according to the orientation of $D(c)$ and calculate how many times $J(N)$ is twisted w.r.t. $B$. Since $D(c)$ is totally real, near each angle of $D(c)$ (corresponding to the endpoints of $d\cup c'$) the field $J(N)$ rotates of $-\frac{1}{4}$ full twist w.r.t. $B$; moreover there is no twist along the vertical edges. Along $c'$, it can be checked that $J(N)$ rotates of $\frac{1}{2}$ full twist w.r.t. $B$ if $c$ is negative and $-\frac{1}{2}$ otherwise. The symmetry is broken on $d$ where it can be checked that $J(N)$ twists w.r.t $B$ always of $-\frac{1}{2}$ full twist. So summing up the various contributions we have that $\nu=-\frac{1}{4}\times 4+\frac{1}{2}-\frac{1}{2}=-1$ at a negative cusp and $\nu=-2$ at a positive one. Hence Theorem \ref{teo:bishop} immediately tells us that $I^++I^-$ vanishes at a negative cusp and it is equal to $-1$ otherwise; moreover it implies that $c_1$ is even at a negative cusp and odd otherwise. It can be checked that the section of the determinant line bundle given by $N\wedge T$ (where $T$ is the field tangent to $\partial D(c)$) performs a total number of twists w.r.t. to $\frac{\partial}{\partial x_1}\wedge \frac{\partial}{\partial x_2}$ contained in $\{ -1,0,1 \}$. Then, if $c$ is negative then $c_1=0$ and hence $I^+=I^-=0$. Finally, by Harlamov-Eliashberg Theorem on annihilation of complex points (see \cite{HE}) we can suppose up to a small perturbation of $int(D(c))$ that $D(c)$ contains no complex points at all. The same argument shows that if $c$ is positive, $D(c)$ can be supposed to contain one hyperbolic point.
 
 } \end{prf} 
\begin{figure}
   \centerline{\includegraphics[width=10.4cm]{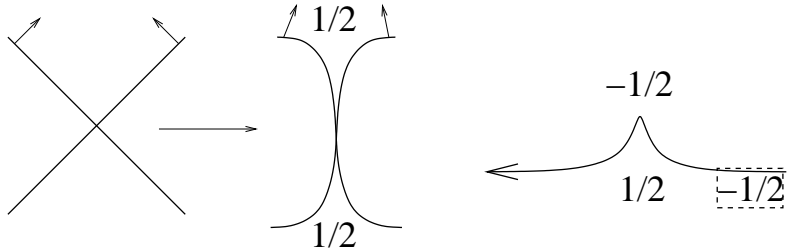}}
  \caption{The Legendrian isotopy on intersection of curves and the local contributions to gleams from cusps. The contribution in the dashed box is assigned to the disc $D(c)$ glued to the cusp to form the mapping cylinder of the projection of the Legendrian curve on $\mr^2$. 
  }\label{fig:polyakdeform}
\end{figure}

We are now able to associate to any oriented Legendrian link $l$ in $\partial M_{\Sigma}$ a branched simple polyhedron $P_l$ embedded in $M_{\Sigma}$ and containing only totally real points and $c^-$ hyperbolic points, where $c^-$ is the number of negative cusps in the projection of the curve. Note that the boundary components of $P_l$ are of two types: those coming from $l$ and those coming from $\partial \Sigma$; in our pictures we draw the former as solid curves and the latter as dashed curves. Since each component of $l$ is equipped with its Thurston-Bennequin framing, we can define a gleam on all the regions of $P_l$ not containing a dashed boundary component. Indeed, we define a section of the normal bundle of each such region near its boundary by using the Thurston-Bennequin framing and the standard construction of the gleam explained in Figure \ref{fig:divergingdirection}; the obstruction to extend this section to the whole region is the half-integer representing the gleam.
In \cite{Po}, an explicit method has been provided to calculate these gleams, which we now recall suitably adapting it to our case. Modify each crossing the diagram of $l$ in $\Sigma$ as showed in the left part of Figure \ref{fig:polyakdeform}; then, on each region of $P_l$ contained in $\Sigma$ subtract the Euler characteristic of the region to the sum of the local contributions assigned to the region as explained in Figure \ref{fig:polyakdeform}. On the (annular) regions of $P_l$ not contained in $\Sigma$ and hence containing one solid boundary component, count the number $c$ of cusps in this Legendrian knot and assign the gleam $-\frac{c}{2}$. 
\begin{rem}
The signs in our constructions are the opposite as those of \cite{Po} because the orientation of the ambient manifold we choose is different.
\end{rem}
We will apply the above constructions in the case when $\Sigma$ is the oriented surface with boundary obtained from a branched polyhedron as follows. 
\begin{defi}
The {\it abstract base} of a branched polyhedron $P$ is the oriented surface $\Sigma (P)$ constructed by picking for each vertex $v_i$ of $P$ an oriented $2$-disc $D_i$ (the one which in Figure \ref{branching} is ``horizontal") and then connecting these discs through rectangular strips (one per each edge of $P$) and extending the orientations of the discs. 
\end{defi}

\begin{figure}
   \centerline{\includegraphics[width=6.4cm]{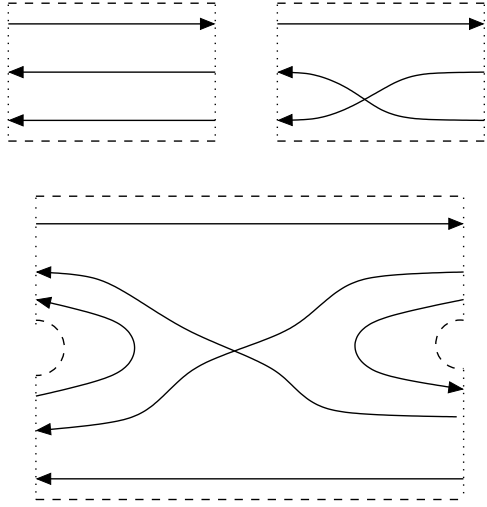}}
  \caption{ In the Figure $\Sigma (P)$ is cut along the dotted curves into discs (the lower block) and strips. The boundary component of the blocks corresponding to $\partial \Sigma (P)$ are dashed, the solid oriented curves are the canonical curves associated to $P$ in $\Sigma (P)$. 
  }\label{fig:basicsteinblocks}
\end{figure}

Note that $\Sigma (P)$ collapses to a $4$-valent graph which is homeomorphic to $Sing(P)$. Starting from the branching of $P$, one can canonically draw a set of oriented curves on $\Sigma (P)$ as explained in Figure \ref{fig:basicsteinblocks} where the upper blocks are contained in the rectangular strips used to construct $\Sigma(P)$ and the lower blocks in the discs $D_i$ and in these last blocks the two intersecting curves correspond to the two regions of $P$ which in the vertex $v_i$ are not ``horizontal" (see Figure \ref{branching}). Let $P''$ be the branched polyhedron obtained by gluing an annulus on each of the above curves along one boundary component and orienting its regions so that the annuli induce the given orientations on the curves and the regions contained in $\Sigma(P)$ are oriented as $\Sigma(P)$. Let moreover $P'$ be the branched polyhedron obtained by puncturing once all the regions of $P$. 
The following is straightforward:
\begin{lemma}\label{lem:p1inp2}
There is a canonical smooth embedding of $P'$ in $P''$. 
\end{lemma}
\subsection{Branched shadows and Stein domains}
In what follows, given a cochain $k$, by $k\leq 0$ we mean that all the coefficients of $k$ are non-positive.
\begin{teo}\label{teo:mainteo}
Let $(P,gl)$ be a branched shadow of a $4$-manifold $M$. If for a non-negative $\mathbb{Z}[\frac{1}{2}]$-lift $UD$ of an ``Up{\rm \&}Down"-cochain the inequality $Eul(P)+gl(P)+\delta UD\leq 0$ holds, then $M$ has a Stein domain structure whose complex structure belongs to the homotopy class encoded by $P$. 
\end{teo}
\begin{cor}\label{cor:branchedspines}
Let $P$ be a branched spine of an oriented $3$-manifold $N$ whose Euler cochain has non-positive coefficients. Then the homotopy class of distribution of oriented $2$-planes on $N$ carried by $P$ contains both a positive and a negative tight contact structure.  
\end{cor}
\begin{rem}
If a $3$-manifold $N$ admits a branched spine satisfying the hypotheses of Corollary \ref{cor:branchedspines}, then $\partial N\neq \emptyset$ and $\chi(\partial N)\leq 0$. Indeed, the maw (see Figure \ref{fig:maw}) can be used to construct a tangent vector field on $\partial N$ containing no singularities with positive index. Hence, in particular, $\partial (N\times [-1,1])$ is a connected $3$-manifold diffeomorphic to the double of $N$.
\end{rem}
\begin{rem}
The above statement on branched spines was conjectured by Benedetti and Petronio in \cite{BP2}. A partial proof of it is possible through the theory of normal surfaces: it is indeed possible to prove that the Chern class of the oriented $2$-planes carried by $P$ is an adjunction class. This is of course a consequence of Mrowka and Kronheimer's result on Seiberg-Witten invariants for $4$-manifolds with boundary (see \cite{KM}).  
\end{rem}
\begin{prf}{3}{
Equip $P$ with $0$ gleam and with the $0$ ``Up\&Down"-cochain. Then, applying Theorem \ref{teo:mainteo} to $P$, one sees that $N\times [-1,1]$ admits a Stein domain structure carried by $P$. Pushing $P$ in $N\times \{1\}$ ($N\times \{-1\}$), we obtain an embedding of $P$ in $N$ as a branched spine and its regular neighborhood is equipped with a positive (negative) tight contact structure inherited from the Stein domain.
}\end{prf}
\begin{prf}{4}{
Let $P'$ be the branched polyhedron obtained by puncturing $P$ along all its regions. Let us fix the choice of the ``Up" and ``Down" regions near each vertex of $P$ which produces the cochain $ud$ whose lift is the given $UD$; this induces a choice of ``Up" and ``Down" regions near each vertex $v_i$ of $P'$. By Lemma \ref{lem:p1inp2}, we can canonically embed $P'$ in a branched polyhedron $P''$ made of blocks as those of Figure \ref{fig:basicsteinblocks}.
The proof is articulated in four steps:
\begin{enumerate}
\item Equip each solid boundary curve of $P''$ with a coorientation and $M_{\Sigma(P)}$ with a Stein domain structure so that $P''\subset M_{\Sigma (P)}$ and  $\partial P'\subset \partial P'' \subset \partial M_{\Sigma (P)}$ is a Legendrian link.
\item Glue $2$-handles corresponding to the regions of $P$ to $P''$ along $\partial P'$ and apply Theorem \ref{teo:gompf-eliashberg}. Call $Q$ the $2$-polyhedron obtained by capping with discs the components of $\partial P'\subset \partial P''$.
\item Collapse the extra regions of $Q$ to get $P$ and calculate the relations between the gleams of the regions of $P$ and those of the regions of $Q$.
\item Clarify the relations between the ambient complex structure and the homotopy class of almost complex structures carried by $P$.
\end{enumerate}

{\bf Step 1.} 
Each region of $P$ corresponds to one of the solid boundary curves of $P''$, so, in each copy of the lower block of Figure \ref{fig:basicsteinblocks} corresponding to a vertex of $P'$, let us associate to the curve corresponding to the ``Up" region the normal vector pointing to its right and to the curve corresponding to the ``Down" region the normal vector pointing to its left.  Let us then complete the choice of the coorientations around the vertex as shown in Figure \ref{fig:basicsteinblocks2} were, at the left we show the case when the region arriving at the vertex in the upper-left edge is ``Up" and at the right we show the other case. Note that, by construction, on each edge touching the vertex, if two curves are oriented in the same direction, they are cooriented in the opposite way. To extend these coorientations along the blocks corresponding to edges of $P'$, we add a number of blocks as the lower one in Figure \ref{fig:basicsteinblocks2} equal to the value of $UD$ on the edge. Since $UD$ is a $\mathbb{Z}[\frac{1}{2}]$-lift of the ``Up\& Down" cochain $ud$ then the coorientations of the curves at the endpoints of each edge will match.

\begin{figure}
   \centerline{\includegraphics[width=11.4cm]{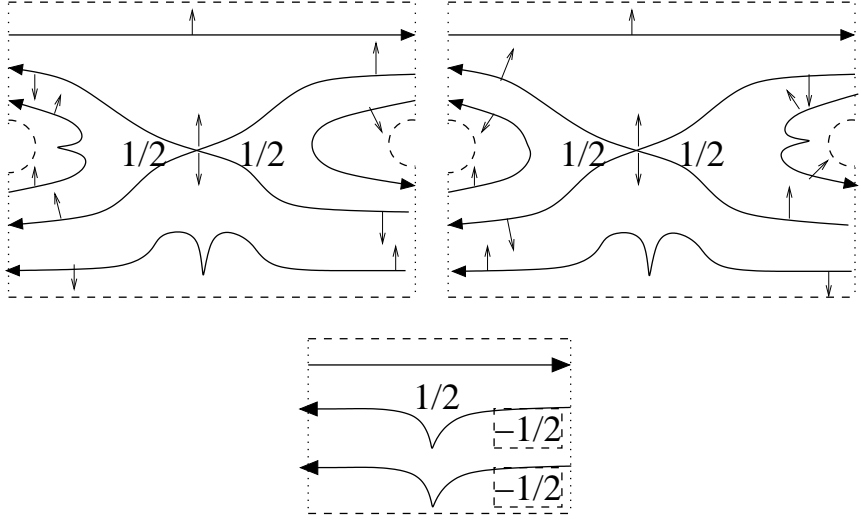}}
  \caption{In the upper part of the figure, we show the coorientations we fix on the blocks corresponding to vertices. In the lower part we show the block we use to contemporaneously revert the coorientations of the two curves oriented the same way on a block corresponding to an edge. The contributions inside the rectangular dashed boxes are to be added to the gleam of the regions of $P''$ containing the cusps. 
  }\label{fig:basicsteinblocks2}
\end{figure}

At the end of the above procedure, we are left with a choice of a coorientation on all the solid boundary curves of $P''$ such that on each intersection point, the Legendrian isotopy of Figure \ref{fig:polyakdeform} produces two cusps whose tangent direction is ``horizontal", as in the internal part  of the upper blocks of Figure \ref{fig:basicsteinblocks2}.

Now that the solid boundary curves of $P''$ are cooriented, they describe a Legendrian link in $\partial M_{\Sigma (P)}$ and the embedded branched polyhedron associated to this link as explained in Subsection \ref{legendriancurves} is exactly $P''$. Note that $M_{\Sigma(P)}$ is a union of $0$ and $1$-handles.

Each region of $P''$ is equipped with gleams calculated by summing the local contributions as explained in Subsection \ref{legendriancurves} and subtracting the Euler characteristic of the region. The regions of $P''$ are either annuli or discs, moreover, since of our choice of the coorientations, for each disc the number of cusps arising from the intersection of two curves on a block corresponding to an edge of $P$ is exactly $2$ so that their total contribution to the gleam of the disc is $1=\chi(Disc)$. Then, the gleam of a region, is the sum only of the contributions given by the blocks of Figure \ref{fig:basicsteinblocks2}. Each region containing a boundary solid curve has gleam equal to $-\frac{c_i}{2}$ where $c_i$ is the number of cusps on that curve.

 {\bf Step 2.} If we glue $2$-handles to the solid boundary components of $P''$ with framing $k_i, i=1,...,n$ w.r.t. the Thurston-Bennequin framing, and with cores $R''_i$ and apply Theorem \ref{teo:gompf-eliashberg} we obtain a Stein structure on the $4$-dimensional thickening of $Q=P''\cup R''_i, \ i=1,..,n$, provided that $k_i\leq -1 \ \forall i, i=1,..,n$.  

The gleam of a region of $Q$ not containing one $R''_i$, is equal to that of the corresponding region of $P''$. The gleam of the region containing $R''_i$ is $k_i-\frac{c_i}{2}$. 

{\bf Step 3.} Collapse $Q$ onto $P=P'\cup R_i, i=1,..,n$, deleting the regions containing the dashed boundary components and obtain an embedding of $P$ in a Stein domain diffeomorphic to its thickening. During this collapse, some regions of $Q$ merge together since they are no longer separated by the singular set: the gleam of the resulting region is the sum of the gleams of these regions in $Q$. In particular, each region $R_i$ of $P$ is the union of the corresponding region $R''_i$ of $Q$ and some other regions contained in the union of the basic blocks of Figure \ref{fig:basicsteinblocks}. 
Let us call {\it external cusp} a cusp not corresponding to one of the two cusps formed in the upper blocks of Figure \ref{fig:basicsteinblocks2} by the intersection of two curves and not coming from a block as the lower one of the figure. It can be checked that, since one of the regions touched by an external cusp is deleted during the collapse and the remaining region has always $\frac{1}{2}$ local contribution from the cusp, the total contribution given by the cusp after the collapse is zero. Then, let $p_i$ be the sum of all the contributions to the regions merging with $R''_i$ coming from one of the two cusps contained inside the upper blocks shown in Figure \ref{fig:basicsteinblocks2} and $c'_i$ be the sum of the contributions coming from the lower block of the figure (including the contributions in the dashed boxes),
so that $gl(R_i)=k_i+c'_i+p_i$. The condition to get a Stein structure of Step 2 can be rewritten as: $gl(R_i)\leq -1+p_i+c'_i$. 

We claim that $1- p_i=Eul_i$. Indeed, reporting the maw (see Figure \ref{fig:maw}) in the lower block of Figure \ref{fig:basicsteinblocks}, we see that it coincides with the vertical direction (i.e. it is parallel to the dotted lines). Then, on each region $R_i$, it rotates by a negative half twist with respect to the outside normal to $R_i$ on $\partial R_i$ exactly when $R_i$ passes though a vertex as one of the two regions which during the collapse of $Q$ on $P$ acquire the $\frac{1}{2}$ contributions from the internal cusps of the upper blocks of Figure \ref{fig:basicsteinblocks2}. Since $Eul_i$ is defined as the number of singularities in $R_i$ of the field extending the maw, and $R_i$ is a disc, $Eul_i$ is equal to $1$ minus $\frac{1}{2}$ the number of times $R_i$ passes through a vertex this way and so  $Eul_i=1-p_i$. Then the above inequality can be rewritten as $Eul_i+gl(R_i)-c'_i\leq 0$ and we are done since $-c'_i=\delta UD_i$ because on each edge $e_j$, we used the lower block of Figure \ref{fig:basicsteinblocks2} exactly $UD(e_j)$-times.

{\bf Step 4.} $P'$ is constructed by gluing regions $R''_i$ to $P''$ which is totally real since the blocks we used contain only positive cusps. Then by Remark \ref{rem:hyperbpoints}, the region $R_i$ of $P$ contains only totally real points and $-1-k_i$ hyperbolic complex points, which we can choose to be positive (the regions of $P$ being oriented by the branching).
Then the ambient complex structure induces no negative complex points on $P$ and so $I^-(P)=0$ (see Definition \ref{def:ipluseiminus}) and Proposition \ref{prop:imenougualediff} applies. 
}\end{prf}
\begin{rem}\label{rem:manysteinstructures}
As in noted in Step 4 and in Remark \ref{rem:hyperbpoints}, the construction gives a totally real embedding of $P$ except on a finite number of hyperbolic complex points contained in the regions. Performing calculations as in Step 3, one sees that the number of these points on a region $R_i$ is $-gl(R_i)-Eul(R_i)-\delta UD(R_i)$. One can choose the numbers $h_i^+$ and $h_i^-$ of positive and negative ones so that $h_i^++h_i^-=-Eul(R_i)-gl(R_i)-\delta UD(R_i)$. Each of these choices gives a Stein domain structure whose complex structure $J$ possibly belongs to a different homotopy class. We can compare the homotopy class of $J$ with the one carried by $P$ (called $J_P$) and express their ``difference" as an integer $2$-cochain (see \cite{Cobr} for details); it can be checked that this difference is $``J-J_P\ "=\sum_i -h_i^- \hat{R}_i$.
\end{rem}
\begin{cor}\label{cor:manysteinstructures}
Let $(P,gl)$ be a shadow satisfying Theorem \ref{teo:mainteo}. With the same notation as in Remark \ref{rem:manysteinstructures}, for each cohomology class $[h]$ represented by a cochain $\sum_i -h_i^-\hat{R}_i$ with $0\leq h_i^-\leq -Eul(R_i)-gl(R_i)-\delta UD(R_i)$ and non-zero in $H^2(P;\mathbb{Z})$, there exists a Stein structure on $M_P$ whose complex structure belongs to the homotopy class which differs of $[h]$ from the one carried by $P$.
\end{cor}
\subsection{Stein neighborhoods of branched shadows}

In this subsection we drop the hypothesis that $M$ is a $4$-handlebody, and by thickening of an embedded polyhedron, we mean the $4$-manifold obtained by applying Theorem \ref{teo:reconstruction} to the polyhedron equipped with the gleams induced by its embedding. We explore the following question:
\begin{question}
Let $P$ be a branched polyhedron embedded in a complex $4$-manifold $(M,J)$ and $gl$ be the gleams induced on it. Under which conditions can be $P$ be isotoped in $M$ so that it admits a neighborhood diffeomorphic to the thickening of $P$ and over which the complex structure induces a Stein domain structure?
\end{question}
The above question was answered when $P$ is an oriented surface by F. Forstneri\v c in \cite{Fo}:
\begin{teo}[\cite{Fo}]\label{teo:forstneric}
An oriented real surface $S$ in $(M,J)$ is isotopic to one whose tubular neighborhood is Stein iff $I^+(S)\leq 0$ and $I^-(S)\leq 0$.
\end{teo}
The natural analogous of the above condition for branched shadows could be that the indices $I^+$ and $I^-$ of each region of $P$ are non-positive. Unfortunately this is not true. Indeed, since by Lemma \ref{lem:branchedHE} the indices $I^+$ and $I^-$ can be viewed as cochains on $P$ which can be independently changed by any coboundary through isotopies of $P$ in $M$, the above assertion reduces to a nearly algebraic one: when a $2$-cochain is equivalent to one expressed by non-positive coefficients? One can prove that each branched shadow can be modified by applying suitable local modifications called ``moves" (see \cite{Cobr}) so that each $2$-cochain is equivalent to a non-positive one. So, since there are $4$-manifolds not admitting any Stein structure (see Lisca's example \cite{Li}), the first generalization we tried is non correct.
We then look at the behavior of the ambient complex structure near the singular set of $P$. 

\begin{defi}
A branched shadow $P$ embedded in a complex manifold $(M,J)$ is in {\it generic position} if it contains a finite number of complex points no one of which is contained in $Sing(P)$ and in {\it good generic position} if furthermore, in each vertex of $P$ the image through $J$ of the maw (the field of Figure \ref{fig:maw}) is the vector transverse to $P$ and lying in the $3$-ball in which the vertex lies (see Definition \ref{def:flat}).   \end{defi}
Each branched polyhedron embedded in $(M,J)$ can be put  in good generic position through an isotopy. Moreover the following holds:
\begin{lemma}\label{lemma:inducedupedown}
If $P$ is in good generic position the ambient complex structure induces on $P$ the positive and negative $2$-cochains $I^{\pm}$ and a $\mathbb{Z}[\frac{1}{2}]$-lift $U$ of an ``Up{\rm \&}Down"-cochain.
\end{lemma}
\begin{prf}{1}{
The first statement has been already treated in the first sections. On each edge we can count the number of twists the projection of $J(maw)$ on the vertical $2$-plane $V(P)$ makes with respect to the direction of Figure \ref{fig:divergingdirection} while running across the edge. Since $M$ is oriented and $P$ is branched, both the edge and $V(P)$ are oriented and we can assign a half-integer to each edge. }\end{prf}
The following is a partial generalization of Theorem \ref{teo:forstneric}:
\begin{teo}\label{teo:genfo}
Let $P$ be a branched shadow in good generic position in a complex manifold $(M,J)$, $gl$ be the gleam induced by the embedding and $I^{+}$, $I^-$ and $U$ be as in Lemma \ref{lemma:inducedupedown}. If there exist $1$-cochains $b^+$ and $b^-$ with values in $\mathbb{Z}$ such that:
\begin{enumerate}
\item The $1$-cochain $U+b^++b^-$ has non-negative coefficients.
\item The $2$-cochains $I^{\pm} +\delta b^{\pm}$ have non-positive coefficients.
\end{enumerate}
Then (up to applying an isotopy to $P$ in $M$) there exists a Stein neighborhood of $P$ in $M$ diffeomorphic to the thickening of $(P,gl)$.  
\end{teo}
\begin{prf}{1}{
The cochains $b^{\pm}$ count the number of times we apply Lemma \ref{lem:branchedHE} to shift the complex points of $P$ and redistribute them so that $I^+$ and $I^-$ are negative on each region. More precisely, let $e_j$ be an edge of $P$, $R_i$, $R_j$ and $R_k$ be the three regions touching it and suppose that $R_i$ is the preferred one. If $b^+(e_j)=q$, we apply Lemma \ref{lem:branchedHE} and let $R_k$ slide over $q$ positive complex points contained in $R_i$ thus shifting them in $R_j$ and $R_k$ (if $q$ is negative then we shift $R_k$ over positive hyperbolic points). We act analogously for $b^-$, shifting $R_i$ over negative points. At the end of these series of isotopies, the cochains of complex points of $P$ have changed and are equal to $I^{\pm}+\delta b^{\pm}$. We claim that each time we shift $R_k$ over an elliptic (hyperbolic) point $p$ (regardless of its sign) then $U(e_j)$ increases by $1$ ($-1$). Indeed let $f$ be an arc contained in $e_j$ and $f'$ and $D$ be as in Figure \ref{fig:slideoverelliptic}. Equip $D$ with the tangent vector field $N$ shown in the figure; by Theorem \ref{teo:bishop}, since $D$ contains only one elliptic (hyperbolic) positive point, $J(N)$ performs $0$ (resp. $-2$) twists w.r.t. the trivial section of the normal bundle of $D$ in $M$. Moreover the maw (see Figure \ref{fig:maw}) can be extended on $f'$ as the field pointing towards $R_i$; then, since the maw rotates of $-1$ twist along $\partial D$ with respect to $N$ then $J(maw)$ rotates of $1$ full twist w.r.t. $J(N)$. Then, since the maw and $N$ coincide (up to sign) on $f$, shifting $R_k$ along $D$ changes the number of twists $J(maw)$ performs on $e$ with respect to the drawn vertical direction of Figure \ref{fig:slideoverelliptic} by $0+1=1$ (resp. $-2+1=-1$).
\begin{figure}
   \centerline{\includegraphics[width=6.4cm]{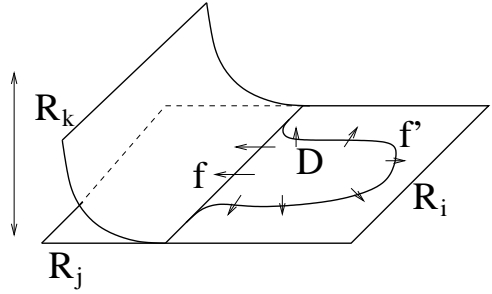}}
  \caption{The vector field $N$ on $D$. 
  }\label{fig:slideoverelliptic}
\end{figure}
After these isotopies, using the maw near $Sing(P)$ we can embed $\Sigma (P)$ in $M$ so that its image is tangent to $P$ along $Sing(P)$ and contains no complex points. Then, on a small neighborhood $M_{\Sigma (P)}$ of $\Sigma (P)$, the ambient complex structure induces a Stein domain structure. We now isotope $P$ near $\Sigma (P)$ keeping it fixed in a small neighborhood  of $Sing(P)$ so that $(P\cup \Sigma (P))\cap M_{\Sigma (P)}$ is the polyhedron $P''$ of Lemma \ref{lem:p1inp2} and $\partial P''\subset M_{\Sigma (P)}$ is a Legendrian link. 
Let us fix an ``Up\& Down" notion near each vertex of $P$ inducing a $1$-cochain whose integer lift is $U+b^++b^-$ (there exists one since $b^+$ and $b^-$ are integer $1$-cochains and $U$ is a lift). Then, near each vertex of $P$ we can isotope $P$ out of a small neighborhood of $Sing(P)$ so that $P\cap \partial M_{\Sigma (P)}$ is Legendrian near the vertices and its projection on $\Sigma (P)$ is as on the upper blocks of in Figure \ref{fig:basicsteinblocks2}. Moreover, one can isotope $P$ near the strips of $\Sigma (P)$ but out of a small neighborhood of $Sing(P)$ so that $P\cap \partial M_{\Sigma (P)}$ is a $3$-uple of curves appearing as in one of the two drawings of Figure \ref{fig:isotopynearedges}, where a non-Legendrian curve $d$ is encoded by a curve $c$ in $\Sigma (P)$ equipped with a non-necessarily orthogonal vector field $v$ in $\Sigma (P)$ and $d$ is given by $c+J(v)$ in local charts around $\Sigma$. Indeed, the idea is to isotope the preferred region by pushing it towards the direction given by $J(maw)$ and push the other two regions in the two opposite directions given by Figure \ref{fig:divergingdirection}; these directions can be encoded on $\Sigma (P)$ by a vector field $v$ such that one region is pushed towards $J(v)$ and the other one towards $J(-v)$.
The number of twists this vector field performs with respect to the maw on each strip composing $\Sigma (P)$ while running on one of the two curves, is equal to the evaluation of $U$ on the edge.
\begin{figure}
   \centerline{\includegraphics[width=8.4cm]{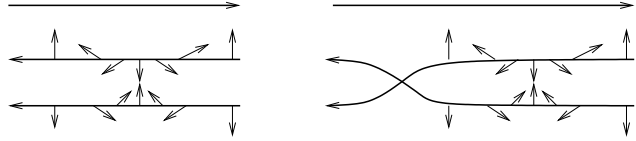}}
  \caption{An example of two possible patterns after the isotopy of $P$ near the edges: in this case the vector field performs $-1$ twist w.r.t. the maw, which in the figure is represented by the vertical direction (it points towards the preferred region). In general any number of twists is possible. 
  }\label{fig:isotopynearedges}
\end{figure}

Then we perform an isotopy of $P$ near $\partial M_{\Sigma (P)}$ which keeps fixed the curves corresponding to the preferred regions of the edges and on the other curves replaces an arc over which the vector field performs $+1$ twist w.r.t. the maw, with a pair of positive cusps as shown in Figure \ref{fig:lagrisotnearedge}.
\begin{figure}
   \centerline{\includegraphics[width=9.4cm]{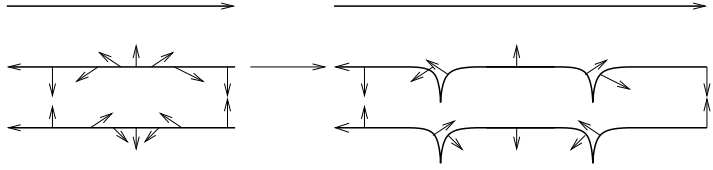}}
  \caption{The isotopy putting in Lagrangian position a curve whose vector field performs $+1$-twist with respect to the maw (which in the figure is vertical). 
  }\label{fig:lagrisotnearedge}
\end{figure}

We end up with a new position of $P$ in $M$ such that $P\cap M_{\Sigma (P)}$ is diffemorphic to $P'$ and its boundary is Legendrian. Moreover, since during the isotopy we did not move $P$ near $Sing(P)$ all the indices of the regions of $P\cap M_{\Sigma (P)}$ in Theorem \ref{teo:bishop} are equal to those of the corresponding regions in the abstract construction of $P'\subset M_{\Sigma (P)}$ and, since we did not use negative cusps, by Lemma \ref{lem:cuspanalisis} $P\cap M_{\Sigma (P)}$ is totally real.

But since $P\cap M_{\Sigma (P)}$ is totally real and we did not move $P$ near $Sing(P)$, all the complex points of $P$ are confined in $P -M_{\Sigma (P)}$ and the indices of the regions of $P$ have not changed during the isotopy and in particular are non-positive. Then, by Forstneri\v c's results (see \cite{Fo}) the distance (using any riemannian metric on $M$) from each component of $P-M_{\Sigma (P)}$ cuts Stein neighborhoods which can be reglued to $M_{\Sigma (P)}$ using the same techniques Eliashberg used in (see \cite{El}). 
}\end{prf}
\subsection{Some applications}
It is easy to find examples of $3$-manifolds admitting both positive and negative Stein-fillable contact structures: consider for instance the $S^1$-bundles with ``low" Euler number on surfaces with sufficiently negative Euler-characteristic.  All these examples have zero Gromov-norm; the following shows that similar examples can be provided also through hyperbolic $3$-manifolds:  
\begin{prop}\label{prop:application1}
For each $n\geq 1$ there exists a closed, oriented, hyperbolic $3$-manifold whose volume is bounded above by $4nVol_{oct}$ (where $Vol_{oct}$ is the volume of the regular ideal hyperbolic octahedron) and admitting $n-1$ positive Stein fillable contact structures and $n-2$ negative ones, all of which have diffeomorphic Stein fillings.
\end{prop}
\begin{prf}{1}{
It is easy to check that for each $n\geq 1$ there exists a branched polyhedron $P$ having $n$ vertices. Up to applying a branched $1\to 2$-move (see \cite{Cobr} for further details on the branched versions of this move) near each vertex of $P$, we can modify it to a branched polyhedron $P'$ having less than $2n$ vertices and such that each region of $P'$ touches an edge for which it is not the preferred region (see Definition \ref{defi:preferred}).
We now claim that we can modify $P'$ in order to obtain a branched polyhedron $P_n$ containing less than $2n$ vertices and only one region.
Indeed, suppose that at an edge $e$ of $P'$ the two non-preferred regions are different, then by exchanging them while passing on $e$ as shown in Figure \ref{fig:hyperbolicstein}, one can modify $P'$ to another branched polyhedron whose number of regions is strictly minor. Hence, along each edge of $Sing(P')$ the two non preferred regions can be supposed to be the same and since each region is the non-preferred region for an edge of $P'$, after repeating the above operation on each edge of $Sing(P')$ we obtain the seeked $P_n$. 

Equipping the region of $P_{n}$ with gleam $0$ or $\frac{1}{2}$ (according to its $\mz_2$-gleam), we can apply Theorem \ref{teo:mainteo} and Corollary \ref{cor:manysteinstructures} and construct Stein domain structures $S^j,\ j=1,\ldots n-1$ on the thickening $M_{P_{n}}$ of $P_{n}$. In order to construct negative contact structures, consider the shadow $\overline{P_{n}}$, obtained by multiplying by $-1$ the gleam of $P_{n}$: $\overline{P_{n}}$ is a shadow of $M_{P_{n}}$ equipped with reversed orientation. Acting as above, we construct Stein structures $\overline{S^j},\ j=1,\ldots n-2$ inducing the seeked negative contact structures on $\partial M_{P_{n}}$. To prove the hyperbolicity of $\partial M_{P_n}$, we use the results of \cite{CT}: $\partial M_{P_n}$ is an integer Dehn surgery of a hyperbolic knot, whose complement has volume bounded above by  $2(2n)Vol_{oct}$ (it is linear in the number of vertices of $P_n$), along meridian curves whose lengths are greater than $6$ (because $P_n$ contains only one regions which touches 6 times each vertex). Then, by Agol and Lackenby's results (\cite{Ag}, \cite{Lac}), $\partial M_{P_n}$ is hyperbolike. To conclude, it is sufficient to observe that $\partial M_{P_n}$ has $b_1\geq 1$ and hence is Haken, whence hyperbolic and to recall that its volume is bounded above by $4nVol_{oct}$ because the volume decreases under Dehn-filling.

\begin{figure}
\psfrag{pr1}{$\partial R_1$}
\psfrag{pr2}{$\partial R_2$}
\psfrag{alfa}{$\alpha$}
\psfrag{pd'}{$\partial D'$}

   \centerline{\includegraphics[width=3.4cm]{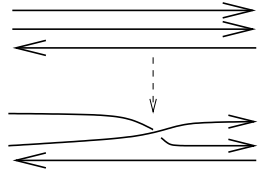}}
  \caption{ The modification ``joining" two distinct non preferred regions along an edge.
  }\label{fig:hyperbolicstein}
\end{figure}
}\end{prf}
In order to provide a last application of Theorem \ref{teo:genfo}, let $P$ a branched $2$-polyhedron $P$ in a $3$ or $4$-manifold $M$; we say that an oriented surface $S$ embedded in $M$ is {\it carried} by $P$ if it is contained in a regular neighborhood of $P$ and its projection on $P$ has always positive differential (i.e. is surjective and orientation preserving). Roughly speaking, such a surface is a union of leaves running parallel to the regions of $P$ and passing over the singular set ``horizontally". When $M$ has dimension $3$, Oertel (\cite{Oe}) provided sufficient conditions on $P$ ensuring that any surface carried by $P$ is incompressible. The following is a $4$-dimensional analogue of that conditions: 
\begin{teo}
Let $(M,J)$ be a complex $4$-manifold collapsing on an embedded branched polyhedron $P$ satisfying the hypotheses of Theorem \ref{teo:genfo} and such that $J$ is homotopic to the almost complex structure carried by $P$ (see Subsection \ref{sub:brshalmostcomplex}). If an oriented and embedded surface $S$ is carried by $P$, then $S$ has the lowest possible genus in its homology class.
\end{teo} 
\begin{prf}{1}{
By Theorem \ref{teo:genfo} a regular neighborhood of $P$ equipped with $J$ is a Stein domain and so its first Chern class $c_1$ satisfies the adjunction inequalities: for each surface $S$ different from a null-homologous sphere, it holds $\chi(S)+[S]^2\leq -|c_1([S])|$. Let $gl(P)$ be the gleam cochain induced by $M$ on $P$, $Eul(P)$ the Euler cochain of $P$; by hypothesis and Corollary \ref{Chern}, $c_1(J)$ is cohomologous to $Eul(P)+gl(P)$.
If $S$ is carried by $P$, up to isotopy we can suppose that $S$ lies very near to $P$ and, pulling-back the maw of $P$ (see Figure \ref{fig:maw}) to $S$, we can construct a tangent vector field on $S$ whose singularities correspond to those of the maw inside the regions of $P$: this shows that $\chi(S)=<Eul(P),[S]>$. Analogously, since the field of vertical $2$-planes transverse to $P$ is also transverse to $S$, we have $[S]^2=<gl(P),[S]>$. Moreover, since $S$ is carried by $P$ we have $[S]\neq 0 \in H_2(P;\mz)$; then the above inequality reads: $$<Eul(P),[S]>+<gl(P),[S]>\leq -|c_1([S])|=-|<Eul(P),[S]>+<gl(P),[S]>|$$
Hence equality holds and $S$ is a minimal genus representative of its homology class.
}\end{prf}
\begin{rem}
The above result is an application of Theorem \ref{teo:genfo}, and applies to branched polyehdra embedded in complex manifolds. A similar result for abstract, non-embedded branched shadows was proved in \cite{Co}: it represents an application of Theorem \ref{teo:mainteo}. More precisely, it was proved that, under suitable combinatorial conditions on an abstract branched shadow $(P,gl)$, a cycle ``carried" by $P$ (i.e. having non-negative coefficients in the canonical presentation of $H_2(P;\mz)$), can be explicitly represented by a minimal genus embedded surface in the thickening $M_{(P,gl)}$ of $(P,gl)$. The proof is based on a construction of ``normal surfaces" in $M_{(P,gl)}$ representing elements of $H_2(M_{(P,gl)};\mz)$.    
\end{rem}

\end{document}